\documentclass[12pt]{amsart}
\usepackage{latexsym,color,amsmath,amsthm,amssymb,amscd,amsfonts}

\newtheorem{thm}{Theorem}

\newtheorem{cor}[thm]{Corollary}

\newtheorem*{remark*}{Remark}

\begin{document}
\title{Tail estimates for sums of variables sampled by a random walk}
\author{Roy Wagner}
\date{October 25, 2006}
\maketitle

\begin{abstract}
We prove tail estimates for variables of the form $\sum_i f(X_i)$, where $(X_i)_i$
is a sequence of states drawn from a reversible Markov chain, or, equivalently, from a random walk on an undirected graph.
The estimates are in terms of the range of the function $f$, its variance, and the spectrum of the graph.
The purpose of our estimates is to determine the number of chain/walk samples which are required
for approximating the expectation of a distribution on vertices of a graph, especially an expander.
The estimates must therefore provide information for fixed number of samples (as in Gillman's \cite{gillman})
rather than just asymptotic information. Our proofs are more elementary than other proofs in the literature,
and our results are sharper. We obtain Bernstein and Bennett-type inequalities, as well as an inequality for
subgaussian variables.
\end{abstract}

\section{Introduction}

One of the basic concerns of sampling theory is economising on the `cost' and quantity of samples required to
estimate the expectation of random variables. Drawing states by implementing a reversible Markov chain or, equivalently, by
conducting a random walk is often considerably
`cheaper' than the standard Monte-Carlo procedure of drawing independent random states. Independence
is indeed lost when sampling by a Markov chain; the empirical average, however, may converge to the
actual average at a comparable rate to the rate of convergence for independent sampling.
This form of sampling is especially useful in the context of random walks on expander graphs.

This approach plays an important role in
statistical physics and in computer science (a concise summary of applications is provided in \cite{hoorylinialwigderson}).
Results concerning the rate of convergence of empirical averages sampled by a random walk, which
hold for a fixed number of samples (rather than just asymptotically), have
been obtained by several authors starting with Gillman's \cite{gillman},
followed by \cite{dinwoodie}, \cite{lezaud} and \cite{leonperron} (for
vector valued functions consult \cite{kargin}). Of these, only \cite{lezaud}
and \cite{kargin} allowed the variance to play a role in their estimates, as is the case in this paper.

This paper is a further step in this direction. We improve known Bernstein-type inequalities,
and prove a new Bennett-type inequality and a new inequality for subgaussian variables.
Our methods are much more elementary than the ones prevailing in the literature,
as we do not apply Kato's perturbation theory to estimate eigenvalues.

Our results were motivated by applications relating to graphs with large spectral gaps (expanders) and
tails which go far beyond the variance (large deviations), such as the recent \cite{artsteinmilman}.
Accordingly, our results are stated for the reversible discrete setting.
Analogues for the continuous and non-reversible settings
can be derived using the simple reduction techniques presented in sections 3.2 and 3.3 of \cite{lezaud}.

\section{The results}

Let $G$ be a finite undirected, possibly weighted, connected graph with $N$ vertices (random walks on such graphs
can represent any finite irreducible reversible Markov chain). Denote by $s$ the stationary distribution of the random walk on the
graph or the equivalent Markov chain. Let $f$ be a function on the vertices of $G$, normalised to have absolute maximum
$1$ and mean $0$ relative to the stationary distribution, namely $\sum_i f(i)s(i) = 0$.
Let $V=\sum_i f^2(i)s(i)$ denote the variance of $f$ with respect to the stationary distribution.
We will think of functions on $G$ as vectors in
$\mathbb{R}^N$ and vice versa, so where $u$ and $v$ are vectors, expressions such as $e^u$ and $uv$ will stand
for coordinatewise operations.

Denote by $P$ the transition matrix of the Markov chain/random walk, such that $P_{ij}$ is the probability of moving from node/state
$j$ to node/state $i$. By the Perron-Frobenius theorem the eigenvalues of this matrix are all real, the top eigenvalue
is $1$ (with $s$ as the only corresponding eigenvector up to scalar multiplication), and the absolute value of
all other eigenvalues is smaller or equal to $1$. Let $\alpha < 1$ be the maximum between the second largest
eigenvalue of $P$ and zero,
and $\beta \leq 1$ the second largest absolute value of an eigenvalue of $P$.

Given a starting distribution $q$, the random variables $X_0,X_1,\ldots$ will denote the trajectory of the random
walk or, equivalently, the states drawn from the Markov chain. $\mathbb{P}_q$ and $\mathbb{E}_q$ will stand for the probability and expectation of events related to this
walk respectively. Let $S_n = \sum_{i=1}^n f(X_i)$. Our concern in this paper is tail estimates for the
distribution of $S_n$.

We will prove inequalities in terms of both $\alpha$ and $\beta$. Note that inequalities in terms of $\alpha$
`cost' an additional multiplicative factor outside the exponent, whereas inequalities in terms of
$\beta$ are useless in the case of $\beta=1$ (i.e. bipartite graphs),
and become relatively poor if $\alpha$ is small and $\beta$ is large, which may be the case.

\begin{thm}\label{maintheorem}
Define
$$\delta(x,r)= x(e^{2r} + V(e^r-1)^2) $$
and
$$\Delta(x,r)= \frac{4xe^{2r}(e^r-1)^2}{1-\delta(x,r)} .$$
In the above setting we get
\begin{align}
\label{betaresult}
\mathbb{P}_q(\frac{1}{n}S_n > \gamma) & \leq
\min_{1 \geq \delta(\beta^2,r), \ r\geq 0}
              \left\|\frac{q}{\sqrt{s}}\right\|_2 e^{-\frac{n}{2}\left[2\gamma r  - V\left(e^{2r} -1 -2r +
                      \Delta(\beta^2,r)\right)\right]}, \\
\nonumber \text{and} & \\
\label{alpharesult}
\mathbb{P}_q(\frac{1}{n}S_n > \gamma) & \leq
\min_{1 \geq \delta(\alpha,r), \ r\geq 0}
              \left\|\frac{q}{\sqrt{s}}\right\|_2 e^{2r}e^{-n\left[2\gamma r  - V\left(e^{2r} -1 -2r +
              \Delta(\alpha,r)\right)\right]} .
\end{align}
\end{thm}

\begin{remark*}
Note that the results are the same up to the factor $\frac{1}{2}$ in the exponent, the multiplicative factor $e^{2r}$
and the replacement of $\beta^2$ by $\alpha$.

Note also that when $\alpha$ goes to $0$, which is effectively almost the case of independent sampling,
$\Delta(\alpha,r)$ also vanishes, and the term we get inside the exponent is the same as the term appearing
in standard proofs of Bennett and Bernstein inequalities for independent variables.
\end{remark*}

The infimum is hard to compute, so we must optimise separately for different parameter regimes.

First we use the above result to derive a Bennett-type inequality (cf. \cite{bennett}).
\begin{cor}\label{bennettresult}
In the above setting,
\begin{align*}
\mathbb{P}_q(\frac{1}{n}S_n > \gamma) & \leq \left\|\frac{q}{\sqrt{s}}\right\|_2
   \exp\left(-\frac{n}{2}tC_{\beta}V\left[(1+\frac{\gamma}{tC_{\beta}V})\log{(1+\frac{\gamma}{tC_{\beta}V})} - \frac{\gamma}{tC_{\beta}V}\right]\right) \\
   & \leq \left\|\frac{q}{\sqrt{s}}\right\|_2 \exp\left(-\frac{n}{2}\gamma\log{\frac{\gamma}{etC_{\beta}V}}\right) ,
\end{align*}
where $t \geq 1$, $C_{\beta}=2\frac{1+\beta^2}{1-\beta^2}$ and provided that $\gamma
\leq(t-1)\frac{1+\beta^2}{\beta^2}V$.
\begin{align*}
\mathbb{P}_q(\frac{1}{n}S_n > \gamma) & \leq (1+\frac{\gamma}{tC_{\alpha}V})\left\|\frac{q}{\sqrt{s}}\right\|_2
   \exp\left(-ntC_{\alpha}V\left[(1+\frac{\gamma}{tC_{\alpha}V})\log{(1+\frac{\gamma}{tC_{\alpha}V})} - \frac{\gamma}{tC_{\alpha}V}\right]\right) \\
   & \leq (1+\frac{\gamma}{tC_{\alpha}V})\left\|\frac{q}{\sqrt{s}}\right\|_2 \exp\left(-n\gamma\log{\frac{\gamma}{etC_{\alpha}V}}\right)
\end{align*}
where $t \geq 1$, $C_{\alpha}=2\frac{1+\alpha}{1-\alpha}$ and provided that $\gamma \leq (t-1)\frac{1+\alpha}{\alpha}V$.
\end{cor}

Our theorem also allows us to reproduce Lezaud's estimates from \cite{lezaud} with improved constants:
\begin{cor}\label{bernsteinresult}
$$\mathbb{P}_q(\frac{1}{n}S_n > \gamma) \leq \left\|\frac{q}{\sqrt{s}}\right\|_2
     e^{-n\frac{1-\beta^2}{1+\beta^2}\frac{\gamma^2}{4(V+\gamma)}} $$
and
$$\mathbb{P}_q(\frac{1}{n}S_n > \gamma) \leq e^{\frac{1-\alpha}{1+\alpha}\frac{\gamma^2}{V+\gamma}}
      \left\|\frac{q}{\sqrt{s}}\right\|_2 e^{-n\frac{1-\alpha}{1+\alpha}\frac{\gamma^2}{2(V+\gamma)}} . $$
\end{cor}

\begin{remark*}
These inequalities imply
$$\mathbb{P}_q(\frac{1}{n}S_n > \gamma) \leq
     e^{\frac{1}{4V}(1-\alpha)\gamma^2}\left\|\frac{q}{\sqrt{s}}\right\|_2 e^{-\frac{n}{8V}(1-\alpha)\gamma^2}$$ for
$\gamma \leq V$, and
$$\mathbb{P}_q(\frac{1}{n}S_n > \gamma) \leq
     e^{\frac{1}{4}(1-\alpha)\gamma}\left\|\frac{q}{\sqrt{s}}\right\|_2 e^{-\frac{n}{8}(1-\alpha)\gamma}$$ for
$\gamma \geq V$. If $\gamma$ is much larger or much smaller than $V$, the constant $8$ can be decreased towards $4$.
Our method allows to improve the constant multiplying $\gamma$ in the denominator, but we will not include the details
because the modification to the proof is cumbersome and straightforward.

The Bennett-type bound improves upon this Bernstein-type result for $\gamma >> V$, provided $\beta$ is small enough.
This allows to see how a smaller $\beta$ reduces the number of required samples.
\end{remark*}

Finally, our technique can be adapted to situations where we have additional information on the distribution of $f$,
such as subgaussian tails. Let $s$ denote here, by abuse of notation, the measure on the vertices of the graph
which corresponds to the stationary distribution.

\begin{thm}\label{subgaussiantheorem}
In the above setting assume also that $s(f \geq t) \leq Ce^{-Kt^2}$ for positive $t$, and
remove the assumption $|f| \leq 1$. Then
$$ \mathbb{P}_q(\frac{1}{n}S_n > \gamma) \leq
\left\|\frac{q}{\sqrt{s}}\right\|_2 e^{-\frac{n}{2}\left(\gamma^2 K- \log{(C\sqrt{\pi K} \gamma + 2)}\right)} ,$$
as long as $\gamma \leq \log{(\frac{1}{2\beta}+\frac{1}{2})}/2K\|f\|_{\infty}$.

We also have
$$ \mathbb{P}_q(\frac{1}{n}S_n > \gamma) \leq
\left\|\frac{q}{\sqrt{s}}\right\|_2 e^{2\gamma K}
     e^{-n\left(\gamma^2 K - \log{(C\sqrt{\pi K} \gamma + 2)}\right)} ,$$
as long as $\gamma \leq \log{(\frac{1}{2\sqrt{\alpha}}+\frac{1}{2})}/2K\|f\|_{\infty}$.
\end{thm}

For some parameter regimes Theorem~\ref{subgaussiantheorem} asymptotically improves upon Theorem~12 from
\cite{artsteinmilman}.

\section{Proofs of results in terms of $\beta$}

In this section we will prove the inequalities involving $\beta$.
Sketches of proofs for inequalities involving $\alpha$ are deferred to the next section.
Before we begin proving we introduce some notation. We will denote $\|u\|_{1/s}=\sum_i \frac{u(i)^2}{s(i)}$
the $\frac{1}{s}$-weighted $\ell_2$ norm on $\mathbb{R}^N$. The inner product associated with this norm
is $\langle u,v \rangle=\sum_i \frac{u(i) v(i)}{s(i)}$. When we refer to the standard $\ell_2$ norm we will use the
notation $\|\cdot\|_2$.

The transition matrix $P$ is not necessarily symmetric, and so its eigenvectors need not be orthogonal (this
would be the case only if $G$ were a regular graph). Reversibility, however, promises that $s_j P_{ij} = s_i
P_{ji}$, and so $P$ is self adjoint and its eigenvectors are mutually orthogonal with respect to
the $\frac{1}{s}$-weighted Euclidean structure. Therefore the $\|\cdot\|_{1/s}$
norm of $P$ restricted to the subspace orthogonal to $s$ is $\beta$, the second largest absolute value of
the eigenvalues of $P$.

\begin{proof}[Proof of Theorem 1.]
The beginning of our proof is identical to that of Gillman's and of those which follow its reasoning.
Take $r\geq 0$. By Markov's inequality
$$ \mathbb{P}_q(\frac{1}{n}S_n > \gamma) \leq e^{-rn\gamma}\mathbb{E}_qe^{rS_n}, $$
where the expectation can be directly expressed and estimated as
\begin{align*}
\sum_{(x_0,\ldots,x_n) \in G^{n+1}} \left( e^{rS_n}q(x_0)\prod_{i=0}^{n-1}(P^T)_{x_i,x_{i+1}} \right)
=      & \ \langle s , (e^{rf}P)^n q \rangle \\
 \leq & \ \|q\|_{1/s} \|Pe^{rf}\|_{1/s}^n .
\end{align*}
Here $e^{rf}$ stands for the diagonal matrix with $e^{rf(i)}$ as diagonal entries, and the inner product is,
we recall, the inner product associated with the $\frac{1}{s}$-weighted $\ell_2$ norm.

At this point Gillman's proof and its variations symmetrise the operator so that its norm will equal its top eigenvalue,
and use Kato's spectral perturbation theory to estimate this eigenvalue.
Our proof, on the other hand, will proceed to simply estimate the norm directly. To do that we will use the equality
$$\|Pe^{rf}\|_{1/s}^2= \max_{\|u\|_{1/s}=1} \langle Pe^{rf}u, Pe^{rf}u \rangle .$$
In order to perform the computation we split the
vector $u$ into stationary and orthogonal components, $u=as +b\rho$,
where $\rho$ is normalised and orthogonal to $s$ in the weighted Euclidean structure. Applying
similar decompositions $e^{rf}s=xs+z\sigma$ and $e^{rf}\rho=ys+w\tau$ we get
\begin{align*}
& \|Pe^{rf}\|_{1/s}^2  =  \\
& \max_{\scriptstyle a^2+b^2=1, \ \rho,\, \sigma,\, \tau}
\left\langle a(xs+zP\sigma) +b(ys+wP\tau), a(xs+zP\sigma) +b(ys+wP\tau) \right\rangle .
\end{align*}
We open the inner product and obtain
\begin{align*}
& \|Pe^{rf}\|_{1/s}^2 = \\
& \max_{\scriptstyle a^2+b^2=1, \ \rho,\, \sigma,\, \tau}
a^2(x^2+z^2\|P\sigma\|^2) + b^2(y^2+w^2\|P\tau\|^2) +2ab(xy+zw\langle P\sigma,P\tau \rangle).
\end{align*}
Denote $p_{\sigma}=\|P\sigma\|^2$, $p_{\tau}=\|P\tau\|^2$ and $p_{\sigma,\tau}=\langle P\sigma,P\tau \rangle$.
Our task is reduced to computing the $\ell_2$ norm of the following $2$ by $2$ symmetric bilinear form:
$$ \left( \begin{array}{cc}
x^2+z^2p_{\sigma}  & xy+zwp_{\sigma,\tau} \\
xy+zwp_{\sigma,\tau} & y^2+w^2p_{\tau}
\end{array}\right) .$$
%
Applying standard computations to derive the norm we get
\begin{align*}
\|Pe^{rf}\|_{1/s}^2 = & \ \frac{1}{2}\biggl[ (x^2+y^2+z^2p_{\sigma}+w^2p_{\tau}) + \\
                      & \ \ \ \ \ \bigl[(x^2+y^2-z^2p_{\sigma}-w^2p_{\tau})^2
                                  + 4z^2w^2(p_{\sigma,\tau}^2 - p_{\sigma}p_{\tau}) \\
                      & \ \ \ \   + 4x^2z^2p_{\sigma} + 4y^2w^2p_{\tau}
                                  + 8xyzwp_{\sigma,\tau}\bigr]^{1/2} \biggr]  \\
                 \leq & \ \frac{1}{2}\biggl[ (x^2+y^2+z^2p_{\sigma}+w^2p_{\tau}) + \\
                      & \ \ \ \ \ \left[(x^2+y^2-z^2p_{\sigma}-w^2p_{\tau})^2
                                  + 4(|xz\sqrt{p_{\sigma}}| + |yw\sqrt{p_{\tau}}|)^2\right]^{1/2} \biggr],
\end{align*}
where we used the Cauchy-Schwarz inequality $p_{\sigma,\tau}^2 \leq p_{\sigma}p_{\tau}$.

To estimate the square root we use the inequality $\sqrt{1+X^2} \leq 1 +\frac{X^2}{2}$. This lead us to
\begin{equation}\label{smallX}
\|Pe^{rf}\|_{1/s}^2
   \leq \ x^2+y^2 + \frac{(|xz\sqrt{p_{\sigma}}| + |yw\sqrt{p_{\tau}}|)^2}{x^2+y^2-z^2p_{\sigma}-w^2p_{\tau}}
\end{equation}
Note that this result depends on assuming that $x^2+y^2 \geq z^2p_{\sigma}+w^2p_{\tau}$.
(For the purposes of the proof of Theorem~\ref{subgaussiantheorem} we require the inequality
\begin{equation}\label{largeX}
\|Pe^{rf}\|_{1/s}^2
   \leq \ x^2+y^2 + |xz\sqrt{p_{\sigma}}| + |yw\sqrt{p_{\tau}}|,
\end{equation}
which is obtained by using $\sqrt{1+X^2} \leq 1 +|X|$, and depends on the same inequality.)

Let us now estimate the components of our formula. We recall that $f$ has mean $0$ with respect to
the stationary distribution $s$ and absolute maximum $1$. We obtain
\begin{align*}
x = \langle e^{rf}s,s \rangle
& = \ 1 + \frac{\langle fs,s \rangle r}{1!}
              + \frac{\langle f^2s,s \rangle r^2}{2!}
              + \frac{\langle f^3s,s \rangle r^3}{3!} + \ldots \\
& \leq \ 1 + V(\frac{r^2}{2!} + \frac{r^3}{3!} + \frac{r^4}{4!} + \ldots)
\leq 1 + V (e^r -1 -r)
\end{align*}
Note also that $|f| \leq 1$ implies that $x \leq e^r$, and that by the arithmeitc-geometric mean
$x = \sum_i s(i) e^{rf(i)} \geq e^{r \sum_i s(i) f(i)} = 1$.

To estimate $y= \langle e^{rf}\rho, s \rangle = \langle e^{rf}s, \rho \rangle$
recall that $\rho$ is normalised and orthogonal to $s$,
and that $\langle fs,\rho \rangle \leq \|fs\|_{1/s} = \sqrt{V}$. We get
\begin{align*}
|y| =  |\langle e^{rf}s, \rho \rangle|
& = \ \langle s,\rho \rangle +
          \frac{\langle fs,\rho \rangle r}{1!} +
          \frac{\langle f^2s,\rho \rangle r^2}{2!} + \ldots \\
& \leq \ \sqrt{V}(r + \frac{r^2}{2!} + \frac{r^3}{3!} + \ldots) \leq \sqrt{V}(e^r-1)
\end{align*}
Note that $x^2+y^2 \leq \|e^{rf}s\|_{1/s}^2 = \langle e^{2rf}s,s \rangle $,
which, as in the computation of $x$ above, is bounded by $1 + V(e^{2r} -1 -2r)$.

Next, using the same estimate as for $y$, we get
$|z| =  |\langle e^{rf}s,\sigma \rangle| \leq \sqrt{V}(e^r-1)$.
For $w = \langle e^{rf}\rho,\tau \rangle$ we use the estimate
$|w|\leq e^r$, which also applies to $x$.

Finally, since the norm of $P$ restricted to the subspace orthogonal to $s$ is $\beta$,
we have $p_{\sigma}, p_{\tau}, p_{\sigma,\tau} \leq \beta^2$

Now we plug our estimates into inequality~(\ref{smallX}), and derive
\begin{align*}
\|Pe^{rf}\|_{1/s}^2
& \leq \ 1 + V(e^{2r} -1 -2r) +
                    \frac{\left( 2e^r \sqrt{V}(e^r-1)\beta\right)^2}
                          {1- \beta^2 e^{2r} - \beta^2 V(e^r-1)^2} \\
& \leq \ 1 + V\left(e^{2r} -1 -2r +
                    \frac{4 \beta^2 e^{2r}(e^r-1)^2}
                         {1- \beta^2(e^{2r} + V(e^r-1)^2)}\right)  \\
& \leq \ \exp\left(V\left(e^{2r} -1 -2r +
                    \frac{4 \beta^2 e^{2r}(e^r-1)^2}
                         {1- \beta^2(e^{2r} + V(e^r-1)^2)}\right)\right),
\end{align*}
as long as $1 \geq \beta^2(e^{2r} + V(e^r-1)^2)$.
To conclude, recall that
$$ \mathbb{P}_q(\frac{1}{n}S_n > \gamma) \leq e^{-n\gamma r}\|q\|_{1/s} \|Pe^{rf}\|_{1/s}^n, $$
so we finally obtain
\begin{equation*}
\mathbb{P}_q(\frac{1}{n}S_n > \gamma) \leq
              \min_{1 \geq \delta(\beta^2,r), \ r \geq 0}
              \left\|\frac{q}{\sqrt{s}}\right\|_2 e^{-n\left[\gamma r  - \frac{1}{2}V\left(e^{2r} -1 -2r +
              \Delta(\beta^2,r)\right)\right]}.
\end{equation*}

\end{proof}

To derive the corollaries and Theorem~\ref{subgaussiantheorem}, we only need to assign suitable values to $r$.
We will restrict to the case $q=s$ in order not to have to carry the $\|\frac{q}{\sqrt{s}}\|_2$ term.

\begin{proof}[Proof of Corollary~\ref{bennettresult}.]
Using the inequalities $(e^r-1)^2 \leq e^{2r}-1-2r$ and
$e^{2r}+V(e^r-1)^2 \leq 2e^{2r}-1$
we bound the expression inside the exponent in inequality~(\ref{betaresult}) by
\begin{align*}
& -\frac{n}{2} \left[2\gamma r  - V(e^{2r} -1 -2r)\left(1 +
    \frac{4\beta^2 e^{2r}}{1- \beta^2 (2e^{2r}-1)}\right)\right] \\
      & \leq -\frac{n}{2}\left[2\gamma r  - V(e^{2r} -1 -2r)
    \frac{1+\beta^2 (2e^{2r}+1)}{1- \beta^2 (2e^{2r}-1)}\right].
\end{align*}

Let $t \geq 1$, $C_{\beta}=2\frac{1+\beta^2}{1-\beta^2}$, $\gamma \leq
(t-1)\frac{1+\beta^2}{\beta^2}V$ and $2r = \log{(1+\frac{\gamma}{tC_{\beta}V})}$. It is easy to verify that
$$\frac{1+\beta^2 (2e^{2r}+1)}{1- \beta^2 (2e^{2r}-1)} \leq t C_{\beta}, $$
so that the above expression will be bounded by $$-\frac{n}{2}\left[2\gamma r - tC_{\beta}V(e^{2r} -1 -2r)\right]. $$

Substituting for $r$ yields
\begin{align*}
& -\frac{n}{2}\left[\gamma\log{(1+\frac{\gamma}{tC_{\beta}V})}  -
               tC_{\beta}V(1+\frac{\gamma}{tC_{\beta}V} -1 - \log{(1+\frac{\gamma}{tC_{\beta}V})})\right] \\
\leq & -\frac{n}{2}tC_{\beta}V
       \left[(1+\frac{\gamma}{tC_{\beta}V})\log{(1+\frac{\gamma}{tC_{\beta}V})} - \frac{\gamma}{tC_{\beta}V}\right] \leq
-\frac{n}{2}\gamma\log{\frac{\gamma}{etC_{\beta}V}}.
\end{align*}
\end{proof}

\begin{proof}[Proof of Corollary~\ref{bernsteinresult}.]
First, we will apply the inequalities $e^{2r} -1 -2r \leq 2r^2e^{2r}$, $e^r-1 \leq re^r$ and $e^{2r}+V(e^r-1)^2 \leq 2e^{2r}-1$
to the exponent in inequality~(\ref{betaresult}). The exponent then turns into
\begin{align*}
-\frac{n}{2}\left[2\gamma r  - Vr^2e^{2r}
              \left(2 + \frac{4\beta^2 e^{2r}}{1- \beta^2 (2e^{2r}-1)} \right)\right]
              & = -n\left[\gamma r  - Vr^2\frac{1+\beta^2}{(1+\beta^2)e^{-2r} - 2\beta^2}\right] \\
              & \leq -n\left[\gamma r  - Vr^2\frac{1+\beta^2}{(1+\beta^2)(1-2r) - 2\beta^2}\right].
\end{align*}
Now we set $r=\frac{1-\beta^2}{1+\beta^2}\frac{\gamma}{2(\gamma+V)}$ and obtain the desired result.
Note that using more careful estimates can lead to a sharper constant multiplying $\gamma$ in
the denominator.
\end{proof}

\begin{proof}[Proof of Theorem 4.]
In this proof we will not assume that $|f| \leq 1$.
For the purpose of this proof we offer a different analysis of the bound $x^2+y^2 \leq \sum_i e^{2rf(i)}s(i)$.
This is simply the expectation of $e^{2rf}$ according to the measure $s$.
We can now evaluate this quantity using the subgaussian information. We get
\begin{align*}
x^2+y^2 & \leq \int_{-\infty}^{\infty} e^{2rt}d\left(-s(f\geq t)\right)
               = \int_{-\infty}^{\infty} 2re^{2rt}s(f\geq t) dt \\
             & \leq 1+ \int_0^{\infty} 2re^{2rt}Ce^{-Kt^2} dt = 1 + C\sqrt{\frac{\pi}{K}} r e^{r^2/K}
\end{align*}

Plugging this estimate into inequality~(\ref{largeX}) together with the simple estimates
$x,w \leq e^{r\|f\|_{\infty}}$ and $y,z \leq e^{r\|f\|_{\infty}}-1$ we obtain
$$\|Pe^{rf}\|_{1/s}^2 \leq 1+ C\sqrt{\frac{\pi}{K}}re^{r^2/K}+2\beta(e^{2r\|f\|_{\infty}}-1) .$$
As noted, inequality~(\ref{largeX}) depends on taking $x^2+y^2 \geq z^2\beta^2+w^2\beta^2$,
which is guaranteed as long as $\beta^2(2e^{2r\|f\|_{\infty}}-1) \leq 1$. We will make the stronger
assumption $\beta(2e^{2r\|f\|_{\infty}}-1) \leq 1$, and obtain the bound
$$\|Pe^{rf}\|_{1/s}^2 \leq (C\sqrt{\frac{\pi}{K}}r+2)e^{r^2/K} . $$

Recalling that
$$ \mathbb{P}_q(\frac{1}{n}S_n > \gamma) \leq e^{-n\gamma r}\|q\|_{1/s} \|Pe^{rf}\|_{1/s}^n, $$
and setting $r=\gamma K$, we conclude the required
\begin{align*}
\mathbb{P}_q(\frac{1}{n}S_n > \gamma) & \leq \|q\|_{1/s} e^{-n\frac{\gamma^2 K}{2}}
                                         (C\sqrt{\pi K} \gamma + 2)^{n/2} \\
                             & = \left\|\frac{q}{\sqrt{s}}\right\|_2 e^{-\frac{n}{2}\left(\gamma^2 K -
                                          \log{(C\sqrt{\pi K} \gamma + 2)}\right)} .
\end{align*}
The condition $\beta(2e^{2r\|f\|_{\infty}}-1) \leq 1$ now reduces to
$\gamma \leq \log{(\frac{1}{2\beta}+\frac{1}{2})}/2K\|f\|_{\infty}$.
\end{proof}

\begin{remark*}
Note that our method allows to increase $\gamma$ as far as
$\log{(\frac{1}{2\beta^2}+\frac{1}{2})}/2K\|f\|_{\infty}$,
where our estimate becomes trivial.
\end{remark*}

\section{Proofs of results in terms of $\alpha$}

The differences between the proofs of results in terms of $\beta$ and $\alpha$ are mostly
computational, so I will only sketch the relevant differences.

\begin{proof}[Proof of Theorem 1.]
As above, our task is to estimate $\|(Pe^{rf})^n\|_{1/s}$. We will use the simple identity
$Pe^{rf} = e^{-\frac{1}{2}rf}e^{\frac{1}{2}rf}Pe^{\frac{1}{2}rf}e^{\frac{1}{2}rf}$
to obtain $$\|(e^{rf}P)^n\|_{1/s} \leq e^r\|(e^{\frac{1}{2}rf}Pe^{\frac{1}{2}rf})^n\|_{1/s}
                                  \leq e^r\|e^{\frac{1}{2}rf}Pe^{\frac{1}{2}rf}\|_{1/s}^n .$$

Since the operator $e^{\frac{1}{2}rf}Pe^{\frac{1}{2}rf}$ is self-adjoint with respect to the weighted
Euclidean structure, we have
$$ \|e^{\frac{1}{2}rf}Pe^{\frac{1}{2}rf}\|_{1/s}=
\max_{\|u\|_{1/s}=1} \langle e^{\frac{1}{2}rf}Pe^{\frac{1}{2}rf}u, u \rangle =
\max_{\|u\|_{1/s}=1} \langle Pe^{\frac{1}{2}rf}u, e^{\frac{1}{2}rf}u \rangle .$$
Decomposing the vectors as in the $\beta$-case (with $\frac{1}{2}r$ replacing $r$) we get
\begin{align*}
& \|e^{\frac{1}{2}rf}Pe^{\frac{1}{2}rf}\|_{1/s}  =  \\
& \max_{\scriptstyle a^2+b^2=1, \ \rho,\, \sigma,\, \tau}
\left\langle a(xs+zP\sigma) +b(ys+wP\tau), a(xs+z\sigma) +b(ys+w\tau) \right\rangle .
\end{align*}
We open the inner product and obtain
\begin{align*}
& \|e^{\frac{1}{2}rf}Pe^{\frac{1}{2}rf}\|_{1/s} = \\
& \max_{\scriptstyle a^2+b^2=1, \ \rho,\, \sigma,\, \tau}
a^2(x^2+z^2\langle P\sigma, \sigma \rangle) + b^2(y^2+w^2\langle P\tau, \tau \rangle)
+2ab(xy+zw\langle P\sigma,\tau \rangle).
\end{align*}
Our task is reduced to computing the $\ell_2$ norm of the same $2$ by $2$ symmetric bilinear form
as in the $\beta$-case, except that $r$ is replaced by $\frac{1}{2}r$, and the definitions of the $p$'s are
now $p_{\sigma}=\langle P\sigma, \sigma \rangle$, $p_{\tau}=\langle P\tau, \tau \rangle$
and $p_{\sigma,\tau}=\langle P\sigma,\tau \rangle$.

The following identity still holds:
\begin{align*}
\|e^{\frac{1}{2}rf}Pe^{\frac{1}{2}rf}\|_{1/s} =
                      & \ \frac{1}{2}\biggl[ (x^2+y^2+z^2p_{\sigma}+w^2p_{\tau}) + \\
                      & \ \ \ \ \ \bigl[(x^2+y^2-z^2p_{\sigma}-w^2p_{\tau})^2
                                  + 4z^2w^2(p_{\sigma,\tau}^2-p_{\sigma}p_{\tau}) \\
                      & \ \ \ \   + 4x^2z^2p_{\sigma} + 4y^2w^2p_{\tau}
                                  + 8xyzwp_{\sigma,\tau}\bigr]^{1/2} \biggr] .
\end{align*}
This time, however, the treatment of the terms inside the square root is slightly more delicate. Let $\lambda_i$
be the eigenvalues of $P$ in descending order, and let $(\sigma^i)_i$ and $(\tau^i)_i$
be the coordinates of $\sigma$ and $\tau$
respectively in terms of the associated orthonormal system. Define
$$p_{\sigma}^+ = \sum_{1>\lambda_i>0} \lambda_i (\sigma^i)^2 \ \
   \text{and} \ \ p_{\sigma}^- = - \sum_{\lambda_i<0} \lambda_i (\sigma^i)^2 ,$$
and decompose $p_{\tau}$ and $p_{\sigma,\tau}$ analogously.
By Cauchy-Schwarz $|p_{\sigma,\tau}^+| \leq \sqrt{p_{\sigma}^+p_{\tau}^+}$,
and the same goes for the $p^-$'s.

All this yields
\begin{align*}
p_{\sigma,\tau}^2 - p_{\sigma}p_{\tau} =
     & \ (p_{\sigma,\tau}^+ - p_{\sigma,\tau}^-)^2 - (p_{\sigma}^+ - p_{\sigma}^-)(p_{\tau}^+ - p_{\tau}^-) \\
   = & \ ((p_{\sigma,\tau}^+)^2 - p_{\sigma}^+p_{\tau}^+) + ((p_{\sigma,\tau}^-)^2 - p_{\sigma}^-p_{\tau}^-)
       + (p_{\sigma}^+p_{\tau}^- + p_{\sigma}^-p_{\tau}^+ - 2p_{\sigma,\tau}^-p_{\sigma,\tau}^+) \\
\leq & \ (\sqrt{p_{\sigma}^+p_{\tau}^-} + \sqrt{p_{\sigma}^-p_{\tau}^+})^2
\end{align*}
and
\begin{align*}
x^2z^2p_{\sigma} + & y^2w^2p_{\tau} + 2xyzwp_{\sigma,\tau} \\
                 = & \ (x^2z^2p_{\sigma}^+ + y^2w^2p_{\tau}^+ + 2xyzwp_{\sigma,\tau}^+) -
                       (x^2z^2p_{\sigma}^- + y^2w^2p_{\tau}^- + 2xyzwp_{\sigma,\tau}^-) \\
              \leq & (|xz\sqrt{p_{\sigma}^+}| + |yw\sqrt{p_{\tau}^+}|)^2 .
\end{align*}
We now combine the two estimates to get
\begin{align*}
4z^2w^2(p_{\sigma,\tau}^2- & p_{\sigma}p_{\tau}) +
        4x^2z^2p_{\sigma} + 4y^2w^2p_{\tau} + 8xyzwp_{\sigma,\tau} \\
  \leq & \ 4\max{(|xz|^2,|yw|^2,|zw|^2)}\left((\sqrt{p_{\sigma}^+p_{\tau}^-} + \sqrt{p_{\sigma}^-p_{\tau}^+})^2 +
                                      (\sqrt{p_{\sigma}^+} + \sqrt{p_{\tau}^+})^2\right) .
\end{align*}
Since $\lambda_2=\alpha$, all $p^+$'s are bounded by $\alpha$.
Note also that $p_{\sigma}^+ + \alpha p_{\sigma}^- \leq \alpha\|\sigma\|_{1/s}^2=\alpha$,
and the same goes for $\tau$. So, in fact, the above is bounded by the expression
\begin{equation*}
4\max{(|xz|^2,|yw|^2,|zw|^2)}\left(\left(\sqrt{p_{\sigma}^+(1-p_{\tau}^+/\alpha)}
                             + \sqrt{p_{\tau}^+(1-p_{\sigma}^+/\alpha)}\right)^2
                             + \left(\sqrt{p_{\sigma}^+} + \sqrt{p_{\tau}^+}\right)^2\right) .
\end{equation*}
Rearranging terms and using Cauchy-Schwarz we get
\begin{align*}
& \left(\sqrt{p_{\sigma}^+(1-p_{\tau}^+/\alpha)} + \sqrt{p_{\tau}^+(1-p_{\sigma}^+/\alpha)}\right)^2
                                                 + \left(\sqrt{p_{\sigma}^+} + \sqrt{p_{\tau}^+}\right)^2 \\
& \leq \ 2\left(p_{\sigma}^+(1-p_{\tau}^+/\alpha) + p_{\tau}^+
         + \sqrt{\left(p_{\sigma}^+ + p_{\tau}^+(1-p_{\sigma}^+/\alpha)\right)
                 \left(p_{\tau}^+ + p_{\sigma}^+(1-p_{\tau}^+/\alpha)\right)}\right)
         \leq 4\alpha .
\end{align*}

So we finally obtain
\begin{align*}
\|e^{\frac{1}{2}rf}Pe^{\frac{1}{2}rf}\|_{1/s} =
                      & \ \frac{1}{2}\biggl[ (x^2+y^2+z^2p_{\sigma}+w^2p_{\tau}) \\
                      & \ + \left[(x^2+y^2-z^2p_{\sigma}-w^2p_{\tau})^2
                                  + 16\alpha\max{(|xz|^2,|yw|^2,|zw|^2)} \right]^{1/2} \biggr] .
\end{align*}

Using the same estimates as in the $\beta$-case,
replacing $r$ by $\frac{1}{2}r$ in the estimates of $x,y,z$ and $w$,
recalling that $p_{\sigma},p_{\tau} \leq \alpha$, and finally changing the bound variable
$r$ into $2r$ we obtain the desired results.
\end{proof}

The other proofs derive from the remark following Theorem~\ref{maintheorem},
which applies to the proof of Theorem~\ref{subgaussiantheorem} as well.

\bigskip

\noindent

\begin{tabular}{@{} l @{\ \ \ \ \ \ \ \ \ \,} l }
Roy Wagner \\
Computer Science Department \\
Academic College of Tel-Aviv-Jaffa \\
2 Rabenu Yeruham Street, Jaffa 68182, Israel \\
and \\
School of Mathematical Sciences \\
Tel Aviv University \\
Tel Aviv 69978, Israel \\
rwagner@mta.ac.il\\
\end{tabular}

\end{document}